\documentclass{TPmod}
\usepackage{url}
\usepackage{setspace}
\usepackage{scrextend}
\usetikzlibrary{arrows.meta}
\usetikzlibrary{cd}
\usetikzlibrary{shapes.misc}
\usetikzlibrary{decorations.pathmorphing}

\makeatletter
\def\@citestyle{\m@th\upshape\mdseries}
\let\citeform\@firstofone
\def\@cite#1#2{{%
  \@citestyle[\citeform{#1}\if@tempswa, #2\fi]}}
\makeatother

\usepackage{hyperref}
\usepackage[capitalize]{cleveref}
\usepackage{tocloft}

\DeclareMathOperator{\Diff}{Diff}

\DeclareMathOperator{\Symp}{Symp}

\def\bP{\mathbb{P}}

\def\bZ{\mathbb{Z}}
\def\bQ{\mathbb{Q}}

\def\bC{\mathbb{C}}

\def\bR{\mathbb{R}}

\usepackage[T1]{fontenc}
\usepackage{mathabx,graphicx}

\numberwithin{equation}{section}
\renewcommand{\theequation}{\thesection.\Alph{equation}}

\title[Mapping-space Gromov--Witten classes]{Mapping-space Gromov--Witten classes and the Fano variety of lines}
\author{Ivan Smith, Abigail Ward}

\date{v1: October 2026}

\email{is200@cam.ac.uk}

\address{Centre for Mathematical Sciences, University of Cambridge, Wilberforce Road, CB3 0WB, U.K.}
\address{School of Mathematical and Physical Sciences, University of Sheffield, 
Hounsfield Road,
Sheffield S3 7RH.}
\email{abigail.ward1@sheffield.ac.uk}

\begin{abstract} {\sc Abstract:} 
Let $Y$ be a smooth cubic 3-fold. We show that $S^2 \times Y$ admits an infinite sequence of cohomologous symplectic forms which are all deformation equivalent but pairwise non-isotopic; McDuff \cite{McDuff:examples} showed the corresponding result on $(S^2)^2 \times T^2$, but these appear to be the first simply-connected and monotone examples. The results are obtained from McDuff's  bordism-flavoured Gromov--Witten invariant, re-interpreted cohomologically in terms of an evaluation map valued in the space of maps from the 2-sphere to a symplectic target. The application then relies on classical properties of the Fano space of lines on $Y$. We also discuss higher-dimensional generalisations.%, which we also describe.%The methods also apply in higher dimensions, and give new obstructions to lifting (homotopically non-trivial families of) diffeomorphisms to symplectomorphisms.
\end{abstract}

\begin{document}
\maketitle

\setlength{\cftbeforesecskip}{2pt}
\tableofcontents

\section{Introduction}\label{sec:intro}

A basic question in symplectic topology is to understand the moduli space of symplectic forms on a given compact oriented smooth manifold $X$, i.e.
\[
\Omega_X := \left\{ \omega \in \Omega^2(X) \, | \, d\omega = 0, \omega^n > 0 \right\}/ \Diff_0(X)
\]
Moser's stability theorem states that, if $\Omega_X$ is non-empty, then  the natural map $\tau: \Omega_X \to H^2(X;\bR)$ is a local homeomorphism. Two symplectic forms are said to be deformation equivalent if they live in the same path-component of $\Omega_X$ and isotopic if they coincide in $\Omega_X$.

In \cite{McDuff:examples}, McDuff introduced a Gromov--Witten type invariant which studied the homotopy class (rather than homology class) of an evaluation map on a moduli space of rational curves. More specifically, she considered a symplectic manifold $X$ with a projection map $X \to S^2$, and a loop $S^1 \to \mathcal{M}_A(J)$ in a moduli space of rational curves in a homology class $A \in H_2(X;\bZ)$, and studied the Hopf invariant of a resulting evaluation
\[
S^1 \times S^2 \stackrel{\ev}{\longrightarrow} X \to S^2
\]
(which she showed homotopically factored through $S^3$).  She used this to give the first examples\footnote{As far as the authors know, apart from elementary variations (taking products of McDuff's examples with further factors), no other examples of the phenomenon have been given since.} of pairs of symplectic forms which were deformation equivalent but not isotopic, on the manifold $(S^2)^2 \times T^2$.  Here the necessary non-trivial loop of rational curves is inherited from the non-trivial topology of the ambient manifold.

This note revisits, re-interprets and generalises her construction, leading to examples of simply-connected monotone symplectic manifolds which display the same phenomenon. A simple case is the following. Let $Y \subset \mathbb{CP}^4$ denote a smooth cubic 3-fold, and let $X = \bC \bP^1 \times Y$ with its monotone symplectic structure.

\begin{thm}\label{thm:main_cubic}
    There is a sequence of cohomologous symplectic forms $\omega_k$ on $X$ which are pairwise non-isotopic but all lie in a fixed symplectic deformation class. 
\end{thm}

Let $Z^2(X;\bR) \subset \Omega^2(X;\bR)$ denote the subspace of closed forms, and $\Omega \subset Z^2(X;\bR)$ the space of symplectic forms, which comes with its natural projection map $p: \Omega \to H^2(X;\bR)$. It follows that there is a path-component of $\Omega$ which meets a fibre of $p$ in infinitely many distinct connected components.

The forms $\omega_k = (f^k)^*\omega$ are obtained as pullbacks of the standard form by powers of a particular diffeomorphism $f$ which involves `twisting' in a family of rational curves.   Our methods also show:

\begin{thm}\label{thm:Torelli}
    $X$ admits a diffeomorphism which acts trivially on cohomology and which preserves the homotopy class of the standard almost complex structure, but which is not isotopic to a symplectomorphism of $(X,\omega)$. 
\end{thm}

If one drops the almost complex  condition, then attractive examples of non-symplectic `Torelli' diffeomorphisms on smooth hypersurfaces in projective space were given by Randal-Williams in the note \cite{ORW}; in that case, the obstruction to being symplectic was that the maps act on the set of homotopy classes of stable almost complex structures without fixed points. In contrast, the map of Theorem \ref{thm:Torelli} is `fragile' in a manner reminiscent of the fragility in  \cite{Seidel:4d}: it becomes isotopic to a symplectomorphism for a symplectic form $\Omega_K$ deformation equivalent to $\omega$ (with $K\gg0$). We also describe  analogously fragile loops \(f_t\) of diffeomorphisms of $X$ which perform a family version of this twist, and which do not lift to $\pi_1\Symp(X,\omega)$ (but would lift if one replaced $\omega$ by $\Omega_K$ for suitable $K$).

As remarked above, McDuff defined her bordism-type Gromov--Witten invariant, of a target symplectic manifold $X$ equipped with a map to $S^2$, as the Hopf invariant of a map $S^1\times S^2 \to X \to S^2$ where $S^1 \subset \mathcal{M}_A(J)$ lay in a space of stable maps to $X$. From that perspective, one key point in the proof of Theorem \ref{thm:main_cubic} is that the Fano surface of lines $F$ on $Y$ has non-trivial rational first cohomology. One also needs to construct a diffeomorphism which `twists' a suitable family of $J$-curves in a way which prevents the map from being symplectic.  This is obtained from a family of rotations of the $S^2$-factor of $X$, defined by a suitable map from $Y$ to $SU(2)$. The fact that such a map provably `twists' the $J$-curves in a way which obstructs the map from being symplectic appeals to the Abel--Jacobi theorem for the cubic 3-fold, identifying $H^3(Y;\bQ)$ and $H^1(F;\bQ)$.

One can take another viewpoint on McDuff's result, which is that if $\mathcal{E}_{A}(J) \to \mathcal{M}_A(J)$ is the universal curve, the optimal evaluation map is not the map 
\[
\mathcal{E}_{A}(J) \to X
\]
usually considered in Gromov--Witten theory, but rather  
 \[
 \mathcal{M}_A(J) \to \operatorname{Map}(S^2,X)/\!\!/SO(3),
 \] 
 i.e. evaluation to the space of unparametrised maps (the notation indicates the homotopy quotient).  Here her original invariant has a cohomological formulation. From this viewpoint, Theorem \ref{thm:main_cubic} generalises to give interesting constraints on symplectic lifts of families of diffeomorphisms of  $\bC\bP^n \times Y_{n+2}^{2n+1}$ where $Y_{n+2} \subset \bC \bP^{2n+2}$ is the degree $n+2$ hypersurface. The required Abel--Jacobi input is in this case due to Shimada \cite{Shimada}. We treat the more general case throughout the main body of the text.

\medskip

\noindent \emph{\textbf{AI declaration.}} 
The authors looked for generalisations of McDuff's examples \cite{McDuff:examples} in the context of $KO$-theoretic Gromov--Witten invariants in 2023-24, motivated by  \cite{AMS2}. In a conversation with ChatGPT Astra on a potential obstruction to lifting Torelli diffeomorphisms to symplectomorphisms coming from compatibility of classical cohomological $A_{\infty}$-structure with quantum eigenspace decomposition,  as described in \cite[Section 2.3]{Smith:ICM}, the LLM produced a Torelli diffeomorphism of $S^2 \times Y$ (which does not fit the desiderata from \emph{op. cit.}). We recognised that this example fit into a mapping-space viewpoint  on McDuff's bordism-flavoured Hopf  invariant, which Astra had not referenced. This manuscript arose from thinking through the implications of that perspective, in  conversation with Astra and each other. %\footnote{Artificial Intelligence has been used significantly, and some readers will judge the adjective essential to the truth of the statement.} 

\noindent \emph{\textbf{Acknowledgements.}}  I.S. is grateful to Daniel Huybrechts for a beautiful lecture series on cubic hypersurfaces given during his visit to Cambridge in 2020.  I.S. is partially supported by Frontier Research Grant EP/X030660/1 from UKRI (in lieu of an ERC advanced grant). A.W. was partially supported by the same grant and by a Royal Society University Research Fellowship held by Cheuk Yu Mak. This work was supported in part by a grant of access to OpenAI models through the \emph{ChatGPT for Academic Researchers} program.

\section{Mapping-space evaluations}\label{Sec:MSE}

\subsection{Mapping space evaluations for minimal-area curves\label{sec:eval}}

Let $(X,\omega)$ be a monotone closed symplectic manifold, and let $A \in H_2(X;\bZ)$ be a class of minimal positive area which is spherical, in the sense that $A \in \mathrm{im}(\pi_2(X) \to H_2(X;\bZ))$. Let $\mathcal{J}_\omega$ denote the contractible space of almost complex structures taming $\omega$. For each $J \in \mathcal{J}_{\omega}$ on $X$ there is a moduli space $\mathcal{M}_A(J)$ of unparametrised $J$-holomorphic rational curves in class $A$. The minimal area assumption implies that this is compact, and that none of its elements can be multiply covered. It follows from classical theory that for regular $J$ (which are generic) it is a smooth manifold.  It carries a universal bundle
\[
\mathcal{E}_A(J) \to \mathcal{M}_A(J)
\]
where the domain is the moduli space of curves in class $A$ with a single marked point.

In Gromov--Witten theory, one usually considers the evaluation map
\[ \operatorname{ev}:\mathcal{E}_{A}(J)\longrightarrow X. \]
In this note we consider instead the evaluation
\[
\operatorname{Ev}:\mathcal{M}_A(J)\longrightarrow \operatorname{Map}_A(S^2,X)\mathbin{/\mkern-6mu/}SO(3)
\]
into the mapping space (the subscript indicates the subspace of maps $v$ for which $v_*[S^2]=A$, and the  quotient indicates the homotopy quotient by the reparametrisation action of the domain). There is a corresponding oriented bordism class
\begin{equation} \label{eqn:the bordism invariant} \operatorname{Ev}_*(\mathcal{M}_A(J))\in \Omega_d^{SO}(\operatorname{Map}_A(S^2,X)\mathbin{/\mkern-6mu/}SO(3)), \end{equation}

where \(d=\dim_{\bR}\mathcal{M}_A(J) = \mathrm{dim}_{\bR}(X) + 2\langle c_1(X),A\rangle - 6 \). The minimal-area condition for $A$ makes this well defined, in the sense that it is independent of $J \in \mathcal{J}_\omega$. 

The mapping space $\operatorname{Map}_A(S^2,X)\mathbin{/\mkern-6mu/}SO(3)$ is functorial under continuous maps of $X$ which preserve the class $A$, hence so is the bordism group. 

The definition of $\mathcal{M}_A(J)$ makes no direct reference to the symplectic structure $\omega$, but only to $J$; the fact that $J$ is taming is relevant only for its compactness.  Consequently, if $f: X \to X$ is a diffeomorphism which preserves the class $A$, then by considering the holomorphic curves for the pull-back almost complex structure $f^*J$, which tames
the form $f^*\omega$ and therefore belongs to $\mathcal J_{f^*\omega}$,  one sees that 
\begin{equation}\label{eqn:pushforward class}
\operatorname{Ev}_*(\mathcal{M}_A(f^*J))
=
(f^{-1})_*\operatorname{Ev}_*(\mathcal{M}_A(J)).
\end{equation}

\begin{lem}\label{lem:obstruct torelli}
    Let $f \in \mathrm{Diff}(X)$ preserve the class $A$. Then if $f$ is homotopic to a symplectic diffeomorphism, the action of $f$ on $\Omega_d^{SO}(\operatorname{Map}_A(S^2,X)\mathbin{/\mkern-6mu/}SO(3))$ preserves the class $\mathrm{Ev}_*(\mathcal{M}_A(J))$.
\end{lem}

\begin{proof}
    Suppose $f$ is homotopic to a symplectomorphism $\phi$. The almost complex structure $\phi^*J$ tames $\phi^*\omega = \omega$. Pick a generic path $\{J_t\}$ of almost complex structures interpolating $J$ and $\phi^*J$ and all taming the fixed form $\omega$. The family moduli space $\mathcal{M}_{A}(J_t)$ then defines a bordism between its ends, i.e. the moduli spaces for $J$ and $\phi^*J$. This shows that $(\phi^{-1})_*$ preserves the class $\mathrm{Ev}_*(\mathcal{M}_A(J))$, and hence so does $\phi_*$. The homotopy between $f$ and $\phi$ shows that the two maps
induce the same action on this bordism group.
\end{proof}

Lemma \ref{lem:obstruct torelli} in particular gives an obstruction to a `Torelli diffeomorphism' of $X$, i.e.\ one acting trivially
on $H^*(X;\mathbb Z)$, being homotopic to a symplectomorphism.
\subsection{McDuff's invariant via mapping spaces\label{sec:hopf}}
%We review McDuff's construction. 

Suppose $\pi:E\longrightarrow B$ 
is an oriented \(S^2\)-bundle and \(s:E\to \bC\bP^j\) restricts to a degree-zero map on every fibre. There is a natural evaluation map
\[
\operatorname{Ev}: B \to \operatorname{Map}_0(S^2,\bC\bP^j)\mathbin{/\mkern-6mu/}SO(3)%\mathcal M_{0}(\bC\bP^j).
\]
where the subscript indicates that we consider maps of degree zero and the quotient is the homotopy quotient. 
The pair $(E,s)$ determines a classifying map
$\operatorname{Ev}$, well defined up to homotopy,
fitting into the commutative diagram
\[
\xymatrix@C=3pc@R=2.5pc{
B \ar[d]_-{\operatorname{Ev}}
&
E \ar[l]_-{\pi}
  \ar[d]^-{\widetilde{\operatorname{Ev}}}
  \ar[r]^-{s}
&
\bC\bP^j \ar@{=}[d]
\\
\operatorname{Map}_0(S^2,\bC\bP^j)
  \mathbin{/\mkern-6mu/}SO(3)
&
\mathcal U_j
  \ar[l]^-{\pi_{j}}
  \ar[r]_-{\operatorname{ev}}
&
\bC\bP^j
}
\]
where $\pi_{j}$ is the universal sphere bundle,
$\operatorname{ev}$ is its evaluation map, and the left
square is a pullback of bundles.

We will understand the cohomology of mapping spaces using their Sullivan minimal models.  In general, given a simply-connected space $X$ with finite rational type, and a minimal model $(\Lambda V,d)$ for $X$, the space $\operatorname{Map}_{0}(S^2, X)$ of null-homotopic maps has a model $(\Lambda(V\oplus\overline{V}),d)$
where $\overline{V}=s^{-2}V^{\geq3}$ and the differential $d$ is extended across the two components. 
This construction is compatible with evaluation; see
\cite[Section~2, Theorem~3.1 and Corollaries~4.7--4.8]{BuijsMurillo}.

\begin{example} \label{ex:sphere} The sphere $S^{2j+1}$ has model
$(\Lambda(\sigma_{2j+1}),0)$, where $\sigma$ represents
the  orientation class. The mapping space
$\operatorname{Map}(S^2,S^{2j+1})$ has model
\[ (\Lambda(\sigma_{2j+1},v_{2j-1}),0),
\qquad v=\overline{\sigma}. \]
Compatibility with evaluation identifies the cohomology
class represented by $v$ with
$\pi_!\operatorname{ev}^*\sigma,$
where $\pi$ and $\operatorname{ev}$ are the canonical projection
and evaluation maps.
The same construction on the universal sphere bundle
defines a class
\[v\in H^{2j-1}\!\left(
\operatorname{Map}(S^2,S^{2j+1})
\mathbin{/\mkern-6mu/}SO(3);\bQ\right),\]
whose restriction is represented by the model generator $v$,
by naturality of fibre integration. This class generates the rational cohomology of the quotient in this degree. \end{example}

 \begin{lem}\label{lem:higher_hopf}
     There is a rational class
\begin{equation}
    u_j\in H^{2j-1}\!\left( \operatorname{Map}_0(S^2,\mathbb{CP}^j) \mathbin{/\mkern-6mu/}SO(3);\mathbb Q \right) \label{higher_hopf}
   \end{equation}
whose restriction to the based mapping space detects the adjoint of the generator of
\(\pi_{2j+1}(\mathbb{CP}^j)\).
 \end{lem}

\begin{proof}
     The Sullivan model of \(\mathbb{CP}^j\) is
\[ (\Lambda(x_2,y_{2j+1}),\;dy=x^{j+1}), \]
so the discussion above implies that $\operatorname{Map}_0(S^2,\mathbb{CP}^j)$ has rational model
\begin{equation} (\Lambda(x_2,y_{2j+1},u_{2j-1}), \;dy=x^{j+1},\;du=0).  \label{eqn:minimal_model} \end{equation}
Here \(u=\overline{y}\) is the first odd cohomology class. On restricting to based maps, the model becomes
$(\Lambda(u_{2j-1}),0)$, so $u$ detects the generator of
$\pi_{2j-1}(\Omega^2_0\bC\bP^j)$, which corresponds under
adjunction to the Hopf generator of
$\pi_{2j+1}(\bC\bP^j)$. The class $u$ survives passage to the \(SO(3)\)-homotopy quotient. In the Serre spectral
sequence of the Borel fibration, a differential out of $u\in E^{0,2j-1}$ lands in
some $E^{4k,\,2j-4k}$, which is spanned by $p_1^k\,x^{j-2k}$. But these classes restrict
non-trivially along the inclusion of constant maps $ 
  \mathbb{CP}^j\times BSO(3)\longrightarrow
  \mathrm{Map}_0(S^2,\mathbb{CP}^j)\mathbin{/\!\!/}SO(3),
$
for which the spectral sequence degenerates and $u$ restricts to zero. The Serre spectral sequence also gives uniqueness of \eqref{higher_hopf}, up to normalisation.
\end{proof}

   % We have worked rationally since that suffices for our purposes.

\begin{defn}
    We define 
    \begin{equation} 
\mathfrak h(E,s)\in H^{2j-1}(B;\mathbb Q) \label{higher_hopf_invt}
\end{equation}
to be the pullback of the class $u_j$ under the classifying map $\mathrm{Ev}$ of $(E,s)$.
\end{defn}

\begin{lem}\label{lem:trivial_if_factors}
    This construction is independent of parametrisations, is invariant under homotopies of \(s\), and extends over bordisms of such families. It vanishes when \(s\) factors through \(B\).
\end{lem}

\begin{proof} Clear.\end{proof}

\begin{ex} \label{ex:mcduff}
   Take $j=1$, and a map $s:E\to\mathbb{CP}^1$.
Then, using that constant maps give a section of the
evaluation fibration,
\[ \pi_1\bigl(\operatorname{Map}_0(S^2,S^2)\bigr)
\cong \pi_3(S^2)\cong\mathbb Z \] and the corresponding generator of \(H^1(\operatorname{Map}_0(S^2,S^2);\mathbb Q)\) extends uniquely to the homotopy quotient by $SO(3)$. 
The natural class
\begin{equation} 
\mathfrak h(E,s)\in H^1(B;\mathbb Q) \label{hopf_invt}
\end{equation}
measures the Hopf invariant around loops in \(B\). This gives a cohomological reformulation of the homotopical invariant studied by McDuff in \cite{McDuff:examples}. 
Indeed, fix a map  $G:S^2\times S^1\longrightarrow S^3 \cong SU(2)$. 
For fixed \(z\in S^2\), we take $E  = S^2\times S^1$ with the map $s: E \to S^2$ given by \(s(x,t)=G(x,t)\cdot z\) 
which (factoring through the Hopf fibration $S^3 \to S^2$) indeed has trivial degree on each $S^2$-fibre. Let $\sigma$ generate $H^3(S^3; \bZ)$. Then with compatible normalisations,
\begin{equation}
    \left\langle\mathfrak h(E,s),[S^1]\right\rangle = \left\langle G^*\sigma,[S^2\times S^1]\right\rangle. \label{eqn:cohomological_hopf} 
\end{equation}
\end{ex}

\begin{ex}
   In the vein of Example \ref{ex:mcduff}, writing $
p_j:S^{2j+1}\longrightarrow\mathbb{CP}^j$ 
for the Hopf projection, Lemma \ref{lem:higher_hopf} replaces  McDuff's examples based on 
\[
S^3\xrightarrow{\deg k}S^3\xrightarrow{p_1}S^2
\]
by examples based on 
\[
S^{2j+1}\xrightarrow{\deg k}S^{2j+1}
\xrightarrow{p_j}\mathbb{CP}^j.
\]
\end{ex}
\medskip

\begin{rmk}
    Our set-up is not well adapted to probe information that is genuinely homotopical even at the level of the mapping space evaluation, and we do not know how to associate evaluation-type invariants to more general elements of the (un)stable stem. 
Suppose $ \alpha:S^m\to S^k$ with $m>k$ 
is a double suspension:
\[
\alpha=\Sigma^2\gamma,\qquad
\gamma:S^{m-2}\to S^{k-2}.
\]
Let $ e:S^{k-2}\to\Omega^2S^k$ 
be adjoint to the identity of \(S^k\), and $
\widehat\alpha=e\circ\gamma$.
The oriented-bordism class
$
[S^{m-2},\widehat\alpha]
\in\widetilde\Omega_{m-2}^{SO}(\Omega^2S^k)$ 
factors through
\begin{equation}
\widetilde\Omega_{m-2}^{SO}(S^{k-2})
\cong\Omega_{m-k}^{SO}(\mathrm{pt}).
\label{nogo}
\end{equation}
The element in \eqref{nogo} is the image of the stable class of \(\gamma\) under
\[
\pi_{m-k}^S\longrightarrow\Omega_{m-k}^{SO}.
\]
That map is zero in positive degree, since a stably framed manifold has vanishing positive-degree Stiefel–Whitney and Pontryagin characteristic classes, and is therefore oriented nullbordant.
\end{rmk}

\subsection{Auxiliary projections}

For general $X$ we have no particular access to $\operatorname{Map}_A(S^2,X)\mathbin{/\mkern-6mu/}SO(3)$. To make progress, we will bridge the discussions of Sections \ref{sec:eval} and \ref{sec:hopf} by making two simplifying assumptions:
\begin{itemize}
    \item we assume that $X$ is equipped with a projection map $\mathrm{pr}: X \to \bC\bP^j$ which has trivial degree on the class $A$, i.e. with the property that $\mathrm{pr}_*(A)=0\in H_2(\bC\bP^j)$; this induces a map
    \[
    \operatorname{Map}_A(S^2,X)\mathbin{/\mkern-6mu/}SO(3) \longrightarrow \operatorname{Map}_0(S^2,\bC\bP^j)\mathbin{/\mkern-6mu/}SO(3)
    \]
    \item we will work rationally, at which point the content of \eqref{eqn:the bordism invariant} is essentially the same as its image in rational homology, which one can detect via pairing it with suitable cohomology classes. On the other hand we can also pull back the class $u_j$ constructed in Lemma \ref{lem:higher_hopf} to give $\mathfrak{h}(\mathcal E_A(J),s) \in H^{2j-1}(\mathcal{M}_A(J);\bQ)$.
\end{itemize}

For simplicity, we will make an even stronger hypothesis, namely that
\[
X = \bC\bP^j \times Y
\]
for a compact smooth Fano variety $Y$ with $c_1(Y) = (j+1)h$, where we write
\[
b \in H^2(\bC\bP^j) \qquad h \in H^2(Y)
\]
for the hyperplane class on $\bC\bP^j$, respectively the integral K\"ahler class on $Y$.  Note that our normalisation makes $X$ monotone, with $c_1(X) = (j+1)(b+h)$, when we take the product symplectic form with $[\omega]=c_1$.  Since $A \in H_2(X;\bZ)$ has zero degree on the projective space factor, we can naturally view $A \in H_2(Y;\bZ)$.

We have the diagram
\begin{equation}\label{eqn:arrows everywhere}
\xymatrix{
\mathcal{E}_A(J) \ar[r]^-{\ev_J} \ar[d]^{\pi} \ar@/^2.0pc/[rr]^{s} & X = \bC\bP^j \times Y \ar[r]^{\mathrm{pr}_{\bC\bP^j}}  & \bC\bP^j  \\
\mathcal{M}_A(J) \ar[r]_-{ \operatorname{Ev}_J} & \operatorname{Map}_A(S^2,X)/\!\!/SO(3)  \ar[r]  & \operatorname{Map}_0(S^2,\bC\bP^j)/\!\!/SO(3)   
}
\end{equation}

An important role will be played by the push-pull operation
\begin{equation} \label{eqn:push-pull} 
\kappa_J = \pi_! \ev_J^* \, : \, H^*(X;\bQ) \longrightarrow H^{*-2}(\mathcal{M}_A(J);\bQ).
\end{equation}
%The K\"unneth formula gives the domain as 
%\[
%H^*(Y) \oplus H^*(Y)\cdot b \oplus \cdots \oplus H^*(Y) \cdot b^j
%\]
Since $\langle b,A\rangle = 0$, we have that $\kappa_J(b)=0$. 
On the other hand, we do obtain an interesting class from $b$ as follows:
\begin{lem}\label{lem:Gysin}
    There is a unique class $\beta_J \in H^2(\mathcal{M}_A(J);\bQ)$ with $\pi^*\beta_J = \ev_J^*(b)$.
\end{lem}
\begin{proof}
    Existence and uniqueness both follow from the Gysin sequence, and the fact that $A$ projects to a class of degree zero, i.e. $\langle b, (\mathrm{pr}_{\bC\bP^j})_*A\rangle = 0$.
\end{proof}

Since $X = \bC\bP^j \times Y$ is a product, for the natural product complex structure $J_0 = J_{\bC\bP^j} \times J_0(Y)$, we have
\[
\mathcal{M}_A(J_0) = \bC\bP^j \times F_A(Y,J_0(Y))
\]
where $F_A(Y)$ is the moduli space of $J_0$-holomorphic rational curves in the class $A$ in $Y$.

We will aim to study invariants of the shape
\[
I_\zeta(J) = \int_{\mathcal{M}_A(J)} \mathfrak{h}(\mathcal{E}_{A}(J),s) \cdot (\beta_J \kappa_J(h^2))^j \cdot \kappa_J(\zeta) \qquad \zeta \in H^*(Y) \subset H^*(X)
\]
(tuning the degree of $\zeta$ appropriately). The precise definition is tailored to the geometry of the examples studied in the next section; in particular, $\kappa_J(h^2)$ arises as a natural polarisation on the variety $F_A(Y,J_0(Y))$.

\section{The main examples}\label{Sec:MainEx}

\subsection{Lines on Fano hypersurfaces}

We continue with the notation of Section \ref{Sec:MSE}, but now specialise to the case where $Y = Y_{j+2} \subset \bC\bP^{2j+2}$ is a smooth  degree $j+2$ hypersurface (symplectically, all such are isomorphic). This has $c_1(Y) = (j+1)h$ where $h$ is the polarisation coming from restricting the ambient hyperplane class.  We will take $A$ to be the class of a line lying on $Y$, which is then obviously minimal area.

\begin{ex}
    If $j=1$ then $X = \bC\bP^1 \times Y_3$ where $Y_3 \subset \bC\bP^4$ is the cubic 3-fold.
\end{ex}

For $X = \bC\bP^j \times Y_{j+2}$, the moduli space of curves in class $A$ has real dimension
\[
\dim_{\bR}(\mathcal{M}_A(J)) = 8j-2.
\]
For simplicity, we will write
\[
F_j = F_A(Y_{j+2}, J_0(Y_{j+2}))
\]
for the moduli space of lines on $Y$, classically known as the `Fano variety of lines'. This is a  complex projective variety, of real dimension $6j-2$. Viewing a line in $Y$ as a line in the ambient projective space $\bC\bP^{2j+2}$ gives an embedding $F_j \hookrightarrow Gr(2,2j+3)$. Let $S$ denote the tautological bundle on the Grassmannian.

\begin{lem}
    Assume that $Y_{j+2}$ is general in the moduli space of hypersurfaces. Then the Fano variety of lines \(F_j\) is smooth of complex dimension
\[ \dim_{\mathbb C}F_j =\dim\operatorname{Gr}(2,2j+3)-(j+3) =3j-1. \]
Moreover, the integrable product complex structure $J_0$ is regular for $A \in H_2(X;\bZ)$.
\end{lem}

\begin{proof}
    The argument applies slightly more generally. Let \(Y\subset\mathbb{CP}^{N}\) be a smooth degree-\(d\) hypersurface cut out by a homogeneous polynomial $P$ and \(u:\mathbb{CP}^{1}\hookrightarrow Y\) have image a line \(L\). For integrable \(J\), the linearised Cauchy–Riemann operator is the Dolbeault operator on \(u^*TY\), so
\[
\operatorname{coker}D_u\cong H^1(L,TY|_L).
\]
The tangent–normal sequence, together with \(H^1(\mathcal O(2))=0\), gives
\[
H^1(L,TY|_L)\cong H^1(L,N_{L/Y}).
\]
Consequently, regularity is exactly the vanishing of this normal-bundle obstruction group.

The normal sequence is
\[
0\longrightarrow N_{L/Y}
\longrightarrow \mathcal O(1)^{\oplus(N-1)}
\xrightarrow{\,dP\,}\mathcal O(d)
\longrightarrow0.
\]
Taking cohomology gives
\[
0\longrightarrow H^0(N_{L/Y})
\longrightarrow H^0(\mathcal O(1))^{\oplus(N-1)}
\longrightarrow H^0(\mathcal O(d))
\longrightarrow H^1(N_{L/Y})
\longrightarrow0.
\]
The middle map is precisely the derivative of the section
\[
s_P\in\Gamma\!\left(\operatorname{Gr}(2,N+1),\operatorname{Sym}^dS^*\right)
\]
whose zero scheme is the Fano scheme \(F(Y)\) of lines in $Y$. Thus
\[
D_u\text{ surjective}
\iff s_P\text{ transverse at }L
\iff F(Y)\text{ smooth of dimension }2N-d-3\text{ at }L.
\]
The last equivalence  follows  from
\[
T_{[L]}F(Y)=H^0(N_{L/Y}),\qquad
h^0(N_{L/Y})-h^1(N_{L/Y})=2N-d-3.
\]

For a degree-\(d\) hypersurface which is general in moduli, the sections \(s_P\) generated by varying the defining polynomial span every fibre of \(\operatorname{Sym}^dS^*\): restriction of homogeneous polynomials to any line is surjective. Generic transversality therefore makes \(s_P\) transverse everywhere. 
\end{proof}

\begin{rmk}
    In the case $j=1$, one does not need to take $Y$ general in moduli. Indeed, specialising the previous argument to the  smooth cubic threefold, \(N_{L/Y}\) has rank two and degree zero. Write
\[
N_{L/Y}\cong\mathcal O(a)\oplus\mathcal O(-a),\qquad a\ge0.
\]
Its inclusion in \(\mathcal O(1)^{\oplus3}\) forces \(a\le1\). Hence the only possibilities are
\[
\mathcal O\oplus\mathcal O,
\qquad
\mathcal O(1)\oplus\mathcal O(-1),
\]
and both have \(H^1=0\). In particular, lines on the cubic of the `second type' \cite{Huybrechts:cubic} are also regular for the integrable complex structure.
\end{rmk}

\subsection{Cylinder maps}

We continue to study the integrable complex structure $J_0 = J_{\bC\bP^j} \times J_0(Y_{j+2})$, with $Y_{j+2}$ general in moduli. Write $\pi: \mathcal{E}_j(Y) \to F_j$ for the universal line restricted to the factor $F_j$ of $\mathcal{M}_A(J_0)$. We have maps
\[
\xymatrix{
\mathcal{E}_j(Y) \ar[r]_{\mathrm{ev}} \ar[d]_{\pi} & Y \\ F_j
}
\]
\begin{lem}
    The class $\ell = \pi_! \ev^*(h^2)$ is ample on $F_j$.
\end{lem}
\begin{proof}
    This is the `Pl\"ucker' polarisation. If $S\to Gr(2,2j+3)$ is the tautological rank two vector
bundle, as above, then the Pl\"ucker polarisation on the
Grassmannian has first Chern class $c_1(S^*)$, and the
projective bundle formula identifies this with the class $\ell$.
\end{proof}

\begin{lem}
The `cylinder map'
\[
T_j = \pi_! \ev^*: H^{2j+1}(Y_{j+2};\bQ) \to H^{2j-1}(F_j;\bQ)
\]
is injective on rational cohomology.
\end{lem}
\begin{proof}
    See Shimada \cite{Shimada}. (When $j=1$ this is a version of the classical Abel--Jacobi theorem for the cubic.)
\end{proof}

\begin{prop}
    \label{lem:hodgeriemann}
    The rational bilinear form  on $H^{2j+1}(Y_{j+2}; \bQ)$ defined by
    \[
    (\eta,\xi)\,  \mapsto \, \int_{F_j}T_j(\eta)\,T_j(\xi)\,\ell^{\,j}
    \]
    is not identically zero.
\end{prop}

\begin{proof}
  Consider the extreme Hodge component
$H^{2j,1}(Y_{j+2})$.
In general, for a degree $d$ hypersurface $W$ in $\bC\bP^{n+1}$, there is an isomorphism
    \[
    H^{n-p+1,p-1}_{prim}(W) = \operatorname{Jac}(P)_{pd-n-2}
    \]
    where $P \in \bC[x_0,\ldots,x_{n+1}] $ is a homogeneous polynomial of degree $d$ cutting out $W$ and  $\operatorname{Jac}(P) = \bC[x_0,\ldots,x_{n+1}] / \langle \partial P / \partial x_j \rangle$ is the Jacobian ring of $P$, graded by setting each $x_i$ in degree 1, cf. \cite{Griffiths}. 
For the hypersurface $Y$ of degree $d=j+2$, and taking $n=2j+1$, we find that $H^{2j,1}(Y) = \operatorname{Jac}(P)_1$ is the degree one  piece, which is necessarily  non-vanishing. Moreover, its injective image under \(T_j\) lies in
$ H^{2j-1,0}(F_j).$ 
For nonzero \(a\) in this image,
\[ \int_{F_j}a\,\overline a\,\ell^j\ne0 \]
by the Hodge--Riemann bilinear relations
\cite[Proposition~3.3.15]{Huybrechts:complex}. This shows the required non-vanishing.
\end{proof}

\begin{rmk}
    In the special case $j=1$, the cylinder map is an isomorphism $H^3(Y;\bQ) \cong H^1(F_1;\bQ)$ and the result follows directly from the Hard Lefschetz theorem on the Fano surface of lines $F_1$.
\end{rmk}

\begin{defn}\label{defn:invariant}
    For $\zeta \in H^{2j+1}(Y_{j+2};\bQ)$ we set
    \begin{equation} \label{eqn:integrand}
I_\zeta(J) = \int_{\mathcal{M}_A(J)} \mathfrak{h}(\mathcal{E}_A(J),s) \cdot (\beta_J \kappa_J(h^2))^j \cdot \kappa_J(\zeta). 
\end{equation}
\end{defn}
Note that the indicated class has degree $(2j-1) + 4j + (2j-1) = 8j-2$ as required.

\begin{lem}\label{lem:independent}
    The characteristic number $I_{\zeta}(J)$ is independent of the choice of $J \in \mathcal{J}_{\omega}$.
\end{lem}

\begin{proof}

Since $A$ has minimal area, a generic path in $\mathcal{J}_\omega$ between two regular choices determines a compact oriented bordism between their moduli spaces. The universal family and evaluation map over this bordism define extensions of all the classes entering the integrand of $I_\zeta$, so the characteristic numbers at the two ends agree.
\end{proof}

\begin{lem}\label{lem:vanishes}
$I_{\zeta}(J_0) = 0$ for every $\zeta \in H^{2j+1}(Y;\bQ)$.
\end{lem}
\begin{proof}
    For the product complex structure, the map $s: \mathcal{E}_{A}(J_0) \to \bC\bP^j$ factors through the base $\mathcal{M}_A(J_0)$ of the universal family, so the result follows from Lemma \ref{lem:trivial_if_factors}.
\end{proof}
Taken together, Lemmas \ref{lem:independent} and \ref{lem:vanishes} look a little disappointing. 
Importantly, the latter need not hold if one evaluates the number for a $J \not\in \mathcal{J}_{\omega}$, for instance a $J$ which tames a different symplectic form.
\begin{lem}\label{lem:obstruction} Let $J\in\mathcal J_\omega$ be regular.
    If $f \in \mathrm{Diff}(X)$ acts trivially on cohomology and $I_{\zeta}(f^*J) \neq 0 $ for some $\zeta$, then $f$ is not isotopic to a symplectic diffeomorphism.
\end{lem}

\begin{proof}
If $f$ were isotopic to a symplectomorphism, Lemma \ref{lem:obstruct torelli} would give
$I_\zeta(f^*J)=I_\zeta(J)$, since $I_\zeta$ is a characteristic
number of the mapping-space bordism class. But
$I_\zeta(J)=I_\zeta(J_0)=0$ by Lemmas \ref{lem:independent} and \ref{lem:vanishes},
a contradiction.
\end{proof}

\subsection{Family twists}

We construct interesting Torelli diffeomorphisms via a `family twist' construction.

\begin{lem}
    Let $\sigma \in H^{2j+1}(S^{2j+1};\Q)$ be our chosen generator. Given any $\eta \in H^{2j+1}(Y; \Q)$, there is some $c\neq 0 \in \bQ$ and a map $\nu: Y_{j+2} \to S^{2j+1}$ with $c\, \eta = \nu^*\sigma  \in H^{2j+1}(Y;\bQ)$.
\end{lem}

\begin{proof}
    Rationally (i.e. after passing to a multiple), every class in $H^{2j+1}(Y)$ may be represented by a map to the sphere, by obstruction theory: the odd-dimensional sphere is rationally Eilenberg--MacLane, so the obstructions are all torsion.
\end{proof}
Replacing $\eta$ by a nonzero multiple if necessary,
we henceforth choose $\nu$ so that $\nu^*\sigma=\eta$. Fix a map $i: S^{2j+1} \to SU(j+1)$ whose composition with projection to the first column has non-zero degree and set $g = i\circ \nu$. For notational convenience, we write $\bar{\imath}\colon S^{2j+1}\to SU(j+1)$ for the pointwise inverse of $i$.

\begin{defn}\label{defn:twist}
     The family twist associated to $\eta$ is
    \[
    f(z,y) := (g(y)\cdot z, y): X = \bC\bP^j \times Y \rightarrow \bC\bP^j \times Y = X.
    \]
\end{defn}

\begin{lem}
    The map $f$ is a diffeomorphism which acts trivially on integral cohomology and preserves the homotopy class of the product almost complex structure.
\end{lem}

    \begin{proof}
         The action map $a:SU(j+1)\times \bC\bP^j\longrightarrow \bC\bP^j$ 
satisfies \(a^*b=\operatorname{pr}_{\bC\bP^j}^*b\), since \(H^2(SU(j+1);\mathbb Z)=0\). Consequently \(f^*b=b\). Since \(f\) covers the identity of \(Y\), the integral K\"unneth decomposition gives
\[ f^*=\mathrm{id}\quad\text{on }H^*(X;\mathbb Z) \]
which verifies the first claim. 
 Relative to
\[ TX=T\mathbb{CP}^{j}\oplus TY, \]
the derivative has the form
\[ df= \begin{pmatrix} A&B\\ 0&I \end{pmatrix}, \]
where \(A\) is complex linear, because each individual rotation is holomorphic. The bundle isomorphisms covering \(f\)
\[ L_t= \begin{pmatrix} A&tB\\ 0&I \end{pmatrix}, \qquad 0\le t\le1, \]
interpolate between a complex-linear isomorphism and \(df\). Pulling back the product complex structure using \(L_t\) gives a homotopy from \(J_0\) to \(f^*J_0\). 
    \end{proof}

\begin{lem}\label{lem:transform}
Under the identification
$\mathcal M_A(f^*J_0)\cong\mathcal M_A(J_0)=\bC\bP^j\times F_j$
given by $C\mapsto f(C)$, we have
\[
\mathfrak h(\mathcal E_A(f^*J_0),s)
=c\,\operatorname{pr}_{F_j}^*T_j(\eta)
\]
for some $c\neq0$.
\end{lem}

\begin{proof}
Under this identification, the universal family is
$\bC\bP^j\times\mathcal E_j(Y)$, and its map to $\bC\bP^j$ becomes
\begin{equation}\label{eqn:action}
s_f(z,q)=\bar{\imath}\bigl(\nu(\operatorname{ev}(q))\bigr)z.
\end{equation}
Let $\operatorname{Ev}_\nu$ classify the family
$\pi:\mathcal E_j(Y)\to F_j$ equipped with
$\nu\circ\operatorname{ev}$. We have a commutative diagram
\[
\xymatrix@C=3pc@R=2.6pc{
F_j
  \ar[d]_-{\operatorname{Ev}_\nu}
&
\mathcal E_j(Y)
  \ar[l]_-{\pi}
  \ar[r]^-{\operatorname{ev}}
  \ar[d]^-{\widetilde{\operatorname{Ev}}_\nu}
&
Y
  \ar[d]^-{\nu}
\\
\operatorname{Map}(S^2,S^{2j+1})
  \mathbin{/\mkern-6mu/}SO(3)
&
\mathcal U_j^{\mathrm{s}}
  \ar[l]^-{\pi_{\mathrm{s}}}
  \ar[r]_-{\operatorname{ev}_{\mathrm{s}}}
&
S^{2j+1}
}
\]
whose left-hand square is a pullback of oriented sphere bundles.
For the class $v$ constructed in Section \ref{sec:hopf},
naturality of fibre integration combined with the commutativity of the diagram gives
\[
\begin{aligned}
\operatorname{Ev}_\nu^*v
&=\operatorname{Ev}_\nu^*
  \bigl((\pi_{\mathrm{s}})_!
         \operatorname{ev}_{\mathrm{s}}^*\sigma\bigr)\\
&=\pi_!\widetilde{\operatorname{Ev}}_\nu^*
        \operatorname{ev}_{\mathrm{s}}^*\sigma\\
&=\pi_!\operatorname{ev}^*\nu^*\sigma\\
&=\pi_!\operatorname{ev}^*\eta
 =T_j(\eta).
\end{aligned}
\]

Let $m:\bC\bP^j\times S^{2j+1}\to\bC\bP^j$ be given by
$m(z,w)=\bar{\imath}(w)z$. The $SO(3)$-equivariant map
\[
\begin{aligned}
\mu_0:\bC\bP^j\times\operatorname{Map}(S^2,S^{2j+1})
&\longrightarrow\operatorname{Map}_0(S^2,\bC\bP^j),\\
\mu_0(z,q)(t)&=m(z,q(t))
\end{aligned}
\]
induces a map $\mu$ on homotopy quotients. By construction,
the twisted family is classified by the map 
$\mu\circ(\mathrm{id}\times\operatorname{Ev}_\nu)$, so
\[
\mathfrak h(\mathcal E_A(f^*J_0),s)
=(\mathrm{id}\times\operatorname{Ev}_\nu)^*\mu^*u_j.
\]
In view of $\operatorname{Ev}_\nu^*v=T_j(\eta)$,
it suffices to show that
\[
\mu^*u_j=c\,\operatorname{pr}_2^*v,
\qquad c\neq0.
\]

From the discussion of Example \ref{ex:sphere} we have 
\[
H^{2j-1}\!\left(
\bC\bP^j\times
\bigl(\operatorname{Map}(S^2,S^{2j+1})
\mathbin{/\mkern-6mu/}SO(3)\bigr);\bQ
\right)
=\bQ\,\operatorname{pr}_2^*v,
\]
so the claim holds for some $c\in\bQ$. To see that $c\neq0$, note that $m([e_1],-)\colon S^{2j+1}\to\bC\bP^j$ factors as
\[
S^{2j+1}
\xrightarrow{\ \operatorname{pr}_1\circ\bar{\imath}\ }
S^{2j+1}
\xrightarrow{\ p_j\ }
\bC\bP^j,
\]
where the first map has non-zero degree by construction, and the second is the Hopf map generating $\pi_{2j+1}(\bC\bP^j)\cong\mathbb{Z}$. Post-composition with $m([e_1],-)$ is the restriction
\[
\mu_0([e_1],-)\colon\Omega^2S^{2j+1}\longrightarrow\Omega^2_0\bC\bP^j,
\]
which by adjunction is therefore non-zero on $\pi_{2j-1}$, sending a generator to a non-zero multiple of the adjoint of $p_j$. Since the restriction of $u_j$ detects that adjoint by Lemma~\ref{lem:higher_hopf}, $\mu^*u_j$ restricts non-trivially to $\{[e_1]\}\times\Omega^2S^{2j+1}$, and $c\neq0$.
\end{proof}

\section{Applications}

\begin{cor}\label{cor:fkdistinct}
    For every \(j\geq1\), the simply-connected monotone Fano manifold \(X_j = \bC\bP^j \times Y_{j+2}\) admits a diffeomorphism \(f\)acting trivially on integral cohomology and preserving the homotopy class of the product integrable complex structure, which is not homotopic to any symplectomorphism. Indeed, the powers of $[f]$ represent pairwise distinct
cosets of the image of
\[
\pi_0\Symp(X,\omega)\longrightarrow\pi_0\Diff(X).
\]
\end{cor}

\begin{proof}
Since $f$ covers the identity on $Y$, it does not change the classes $\kappa_J(\cdot)$ entering into the integrand in \eqref{eqn:integrand}; using Lemma \ref{lem:Gysin} shows that it does not change any class  other than $\mathfrak{h}$. Combining this with Lemma \ref{lem:transform}, we have that for any $\zeta \in H^{2j+1}(Y)$  
\begin{equation} I_{\zeta}(f^*{J_0}) = c\int_{F_j}T_j(\eta)\,T_j(\zeta)\,\ell^{\,j}\label{eqn:high_dim_twisted} \end{equation}
(where integrating out the factors of $\beta_J$ against the $\bC\bP^j$-factor of the moduli space $\mathcal{M}_A(J_0)$ reduces the integral to one on the Fano variety $F_j$). 
By Proposition \ref{lem:hodgeriemann}, we may choose
rational classes $\eta,\zeta$ for which this integral
is nonzero. After replacing $\eta$ by a nonzero integral
multiple, choose $\nu$ as above with $\nu^*\sigma=\eta$.
For the resulting twist $f$, we therefore have
$I_\zeta(f^*J_0)\ne0$; moreover, after iteration we have
\begin{equation} I_\zeta((f^k)^*J_0)= k\, I_\zeta(f^*J_0) \end{equation}
and all of these numbers are distinct. The rest of the proof follows from bordism invariance of the characteristic number and the vanishing result of Lemma \ref{lem:vanishes}.
\end{proof}

\begin{prop}\label{prop:def_equivalent}
    For $k\in \bZ$, the forms \(\omega_k = (f^k)^*\omega\) all represent \((j+1)(b+h)\), and all lie in the same path component of the space of symplectic forms on $X$. 
\end{prop}

\begin{proof}
    Only the second statement requires proof: the argument  works for any fibrewise rotation of \(\bC\bP^j\times Y\). 
    
    For \(K\ge1\), put $\Omega_K=\omega_{\bC\bP^j}+K\omega_Y.$ Since \(f\) covers the identity on \(Y\), $ f^*\omega_Y=\omega_Y$. Consider the closed forms
\begin{equation}
    \Theta_{K,t} =(1-t)\omega_{\bC\bP^j}+t f^*\omega_{\bC\bP^j}+K\omega_Y, \qquad 0\le t\le1. \label{interpolate}\end{equation}

Their restrictions to every \(\bC\bP^j\)-fibre are exactly \(\omega_{\bC\bP^j}\), because the action of $SU(j+1)$ preserves the usual Fubini--Study form. For sufficiently large \(K\), we claim all the forms \eqref{interpolate} are symplectic. Indeed, writing
\[ \alpha_t=(1-t)\omega_{\bC\bP^j}+t f^*\omega_{\bC\bP^j}, \]
the top exterior power on $\bC\bP^j \times Y$ satisfies

\begin{equation} (\alpha_t+K\omega_Y)^{3j+1} = C\cdot K^{2j+1}\omega_{\bC\bP^j}^j\wedge\omega_Y^{2j+1}+O(K^{2j}) \label{positive top power} \end{equation}
for a positive constant $C>0$.
The leading term is a positive volume form, and dominates all other terms for sufficiently large $K$, uniformly over the compact space \(X\times[0,1]\). Thus \eqref{positive top power} is positive for all sufficiently large \(K\).

We now concatenate three paths:
\[ \omega \ \longrightarrow\ \Omega_K \ \longrightarrow\ f^*\Omega_K \ \longrightarrow\ f^*\omega.  \]
The first is the path of product forms obtained by increasing \(K\). The middle is \eqref{interpolate}. The last is the pullback of the reverse of the first path. This proves deformation equivalence, and the same argument applies to every \(f^k\).
\end{proof}

\begin{rmk}
The path \eqref{interpolate} consists of cohomologous symplectic
forms. Moser's theorem therefore gives an isotopy $\psi_t$,
with $\psi_0=\mathrm{id}$ and
\[
\psi_1^*f^*\Omega_K=\Omega_K.
\]
Thus $f\circ\psi_1$ is a symplectomorphism of $(X,\Omega_K)$
isotopic to $f$. The same argument applies parametrically to compact families of
such fibrewise rotations, with $K$ chosen uniformly and the
homotopy fixed at parameters where the rotation is the identity. \end{rmk}

\begin{cor}
    The forms \(\omega_k\) are pairwise non-isotopic.  Each is monotone and Kähler for the pulled-back Fano complex structure \((f^k)^*J_0\); these complex structures are mutually homotopic.
\end{cor}

Recall that here  “non-isotopic” means that there is no path of cohomologous symplectic forms joining distinct members. The forms are manifestly symplectomorphic by diffeomorphisms not lying in  the identity component of the diffeomorphism group.

\begin{prop}
For $X = \bC\bP^j \times Y_{j+2}$ and 
    for each \(r \in \{1,3,\ldots,2j-1\}\), there is an element of infinite order in
\[ \operatorname{coker}\bigl( \pi_r\mathrm{Symp}(X,\omega)\to\pi_r\mathrm{Diff}_0(X) \bigr). \]
\end{prop}

In fact, its image remains nonzero in the cokernel of the comparison with \(\pi_r\mathrm{HtEq}_0(X)\).

\begin{proof}
    We first describe the required family of diffeomorphisms. Let $a = (2j+1-r)/2$. Let \(u\) generate \(H^r(S^r;\mathbb Z)\). Choose a based family
\[ G:S^r\times Y_{j+2}\longrightarrow SU(j+1), \qquad G(s_0,y)=1, \]
whose top odd class pulls back to
$N\,u\cdot h^{a}$ with $N\ne0$ and whose lower odd classes pull back to zero. Such a $G$ can be constructed by first mapping
\[ S^r\wedge (Y_{j+2})_+\longrightarrow S^{2j+1} \]
to realise a suitable nonzero multiple of \(u h^a\), and then composing with a map \(S^{2j+1}\to SU(j+1)\) which is nonzero on the top rational homotopy group. 
Define
\begin{equation} f_s(z,y)=(G(s,y)z,y). \label{last family} \end{equation}
This is a based \(S^r\)-family of Torelli diffeomorphisms. The homotopy which scales away the off-diagonal block of \(df_s\) works simultaneously in \(s\), and is fixed at \(s_0\). Thus the family also admits a coherent homotopy from \(f_s^*J_0\) to \(J_0\).

Denote by \(\mathfrak h_j\) the degree-\((2j-1)\) class obtained by pullback of the class $u_j$ from Lemma \ref{lem:higher_hopf}, with respect to the evaluation twisted by $f_s^{-1}$. Applying the calculation of Lemma 3.15 over $S^r\times F_j$,
we obtain, for some fixed $c_j\ne0$ depending only on our
normalisations,
\[ \mathfrak h_j =c_j\, N\,u\cdot \pi_!\ev^*(h^a),  \]
of degree $ r+(2a-2)=2j-1.$
We now consider a parametrised version of our characteristic class invariant, namely
\begin{equation} \mathcal I_{j,a}(J)= \int_{S^r\times\mathcal{M}_A(J)} \beta^j\,\mathfrak h_j\,\kappa(h^2)^{\,3j-a} \label{eqn:final_char} \end{equation}
where we integrate a form of  total degree  $2j+(2j-1)+2(3j-a) =r+8j-2.$

Using the product family, \eqref{eqn:final_char} becomes
\begin{equation} \mathcal I_{j,a}(f^*J_0) =c_j \, N\int_{F_j} \pi_!\ev^*(h^a)\,\ell^{\,3j-a}. \label{eqn:final_compute} \end{equation}
We claim this integral is nonzero. The class \(\pi_!\ev^*(h^a)\) is the incidence class of lines meeting a general codimension-\(a\) linear subspace. For \(1\le a\le j\), it is represented by a nonempty effective cycle of the expected dimension: \(Y_{j+2}\) is covered by lines, and a general incident line meets the linear subspace at one point. Intersecting with the appropriate power of the ample Plücker class gives a positive degree, proving the claim.

Suppose \(f_s\) were homotopic, as a based family of maps, to a family of symplectomorphisms \(\phi_s\). The inverse families are then homotopic, so their twisted
evaluation maps define the same characteristic number \eqref{eqn:final_char}. Choose a generic extension of the family $J_s=(\phi_s)^*J_0$ to a family of compatible almost complex structures 
 over \(D^{r+1}\). The associated  parametrised moduli space is a compact oriented nullbordism of the family associated to $\phi_s$. All three cohomology classes in \eqref{eqn:final_char} extend over this bordism. Hence \eqref{eqn:final_char} must vanish, contradicting \eqref{eqn:final_compute}. Finally,  the characteristic number is additive under addition in \(\pi_r\), which proves infinite order in the cokernel.
\end{proof}
\begin{rmk}
    The contradiction works even if the homotopy between \(f_s\) and \(\phi_s\) passes through homotopy equivalences rather than diffeomorphisms, justifying the final sentence of the statement of the Proposition.
\end{rmk}

As is well known, such results also have consequences for moduli spaces of symplectic forms. Let
\[ \mathcal S_{[\omega]} = \{\Omega\text{ symplectic on $X_j$}:[\Omega]=[\omega]\}, \qquad \mathcal S=\{\Omega\text{ symplectic}\}, \]

based at \(\omega\), and let $\mathcal{S}_{[\omega]}^0$ denote the component containing $\omega$. The family \eqref{last family} gives
\[ \Phi:S^r\longrightarrow\mathcal S_{[\omega]}, \qquad \Phi(s)=f_s^*\omega.  \]
The Moser fibration
\[ \mathrm{Symp}(X_j,\omega)\cap\mathrm{Diff}_0(X_j) \longrightarrow \mathrm{Diff}_0(X_j) \longrightarrow \mathcal S_{[\omega]}^{\,0} \]
shows the class \([\Phi]\) has infinite order in $\pi_r \mathcal{S}_{[\omega]}$.  Nevertheless, \(\Phi\) is nullhomotopic in \(\mathcal S\), by an argument directly analogous to that given in Proposition \ref{prop:def_equivalent}. Thus, 
\[  \ker\!\left[ \pi_r\mathcal S_{[\omega]} \longrightarrow\pi_r\mathcal S \right] \text{ contains an infinite-order element}  \]
for every odd \(r=1,3,\ldots,2j-1\).

Finally, we obtain an obstruction to the existence of families of symplectic submanifolds, or pseudo-holomorphic projection maps in a particular homotopy class, in the following sense. For simplicity, we just state this in the simplest case $j=1$, since symplectic submersions to $S^2$ arise naturally elsewhere in the study of the Seidel element. We recall the sequence $\omega_k = (f^k)^*\omega$ of symplectic forms constructed above.

\begin{prop}
    Let $X = S^2 \times Y_3$ and let $r:X\to S^2$ be the natural projection. For \(k\ne0\), there cannot exist an \(\omega_k\)-tame almost complex structure \(J\) and a map \(r'\simeq r\), with $r': X \to \bC\bP^1$ $J$-holomorphic.
\end{prop}

\begin{proof}
Suppose we have some $J$ for which $r': X \to \bC\bP^1$ is $J$-holomorphic with $J$ taming $\omega_k$, $k \neq 0.$  The moduli space $\mathcal{M}_A(J)$ is compact but not necessarily smooth. Since $(r')_*(A)=0$, every $J$-curve $v$ in $\mathcal{M}_A(J)$ is mapped to a point by $r'$. Take a sequence of regular  $\omega_k$-tame almost complex structures $J_n$ converging to $J$, and fix a metric on $S^2$ with convexity radius $\rho$. Then there exists some $N>0$ such that for all $n > N$, all $v \in \mathcal{M}_{A}(J_n)$ satisfy 
\[ \mathrm{diam}(r' \circ v) < \rho; \]
if not, we could choose a sequence of $v_n \in \mathcal{M}_{A}(J_n)$ with $r' \circ v_n$ violating this bound, which by Gromov compactness admits a subsequence converging to an element $v_\infty \in \mathcal{M}_{A}(J)$, a contradiction. 

Now fix some $n>N$, and define a map 
\[ c: \mathcal{M}_{A}(J_n) \to S^2 \]
by setting $c(v)$ to be the centre of the smallest geodesic ball containing the image $r' \circ v$. The map $c$ is well defined since the intersection of any two balls with distinct centres but the same radius is contained within a ball of strictly smaller radius. It is also continuous, since the images $r' \circ v$ vary continuously across $\mathcal{M}_{A}(J_n)$ in the Hausdorff topology.

As before let $ \pi: \mathcal{E}_A(J_n) \to \mathcal{M}_{A}(J_n)$ denote the universal sphere bundle. Contracting along geodesic flows defines a homotopy 
$ r' \circ \operatorname{ev} \simeq c \circ \pi,$ 
so by Lemma \ref{lem:trivial_if_factors}, the associated invariant $\mathfrak{h}(\mathcal{E}_A(J_n),  r' \circ \operatorname{ev})$ vanishes and  by homotopy invariance $\mathfrak{h}(\mathcal{E}_A(J_n),  r \circ \operatorname{ev})$ vanishes as well. This implies that $I_\zeta(J_n)=0$ for all $\zeta \in H^3(Y_3; \bQ)$, contradicting the non-vanishing exhibited in the proof of Corollary \ref{cor:fkdistinct}. 
\end{proof}

For the same reason, there is no \emph{symplectic} submersion $r': (X,\omega_k) \to S^2$ in the homotopy class of the projection map $r$ (i.e. no submersion with symplectic fibres; the symplectic orthogonal connection would allow one to choose a compatible \(J\) making \(r'\) holomorphic).
    
\bibliographystyle{amsalpha}
{\footnotesize{\bibliography{bib}}}

\end{document}